\documentclass[12pt]{amsart}
\usepackage{amssymb}
\usepackage{amsmath}
\usepackage[all]{xy}
\usepackage{xspace}
\usepackage{enumerate}
\usepackage{color}

\newtheorem{thms}{Theorem}[section]
\newtheorem{props}[thms]{Proposition}
\newtheorem{lms}[thms]{Lemma}
\newtheorem{defs}[thms]{Definition}
\newtheorem{remar}[thms]{Remark}
\newcommand{\GG}{\mathcal{G}_0}
\newcommand{\GQ}{\mathcal{G}_0\otimes\QQ}
\newcommand{\DD}[1]{\mathbb{I}_{#1}^{\sim}}

\newcommand{\infinv}{\underline{K}_{\bf n}}
\newcommand{\Dilim}{{\displaystyle\lim_{\longrightarrow}}_{\Delta}}
\newcommand{\ppDilim}{{\displaystyle\lim_{\longrightarrow}}_{\Delta'}}
\newcommand{\kmapk}[2]{\kappa_{{#1},{#2}}}

\newcommand{\pKone}[2][n]{\K_1({#2};\ZZp{#1})}
\newcommand{\pKi}[2][n]{\K_i({#2};\ZZp{#1})}
\newcommand{\pKu}[2]{\K_0({#2};\ZZ\oplus\ZZp{#1})}

\newcommand{\pKs}[1]{\K_{0}({#1};\QQ\oplus\QQ/\ZZ)}
\newcommand{\ZZp}[1]{\ZZ/{#1}}

\newcommand{\po}{^+}

\newcommand{\comak}[2]{\chi_{{#1},{#2}}}
\newcommand{\id}{\operatorname{id}}

\newcommand{\infer}{\mbox{$\Longrightarrow$}}
\newcommand{\egmap}[1]{\lfloor{#1}\rfloor}
\newcommand{\semap}[1]{\lceil{#1}\rceil}
\newcommand{\Ideal}{\mathrm{I}}
\newcommand{\tor}{\operatorname{tor}}
\newcommand{\A}{\mathfrak{A}}
\newcommand{\B}{\mathfrak{B}}
\newcommand{\BK}[1][n]{\overline{K}_{#1}}
\newcommand{\BKf}{\overline{K}_{\mathbf n}}

\definecolor{mygray}{gray}{0.7}

\newcommand{\nbox}{\framebox[0.55em]{\raisebox{-0.02em}[0.15em][0em]{$\scriptstyle n$}}}
\newcommand{\fnbox}{{\framebox[0.55em]{\raisebox{-0.02em}[0.15em][0em]{$\scriptstyle
      \mathbf n$}}}}

\newcommand{\gz}{\ensuremath{\mathrm{G}_0}}

\newcommand{\go}{\ensuremath{\mathrm{G}_1}}
\newcommand{\gn}{\ensuremath{\mathrm{G}_n}}

\newcommand{\gx}{\ensuremath{\mathrm{G}_*}}
\newcommand{\gd}{\ensuremath{\mathrm{G}_{\nbox}}}

\newcommand{\kz}{\ensuremath{\mathrm{K}_0}}
\newcommand{\ko}{\ensuremath{\mathrm{K}_1}}

\newcommand{\kx}{\ensuremath{\mathrm{K}_*}}
\newcommand{\kd}{\ensuremath{\mathrm{K}_{\nbox}}}
\newcommand{\hz}{\ensuremath{\mathrm{H}_0}}
\newcommand{\ho}{\ensuremath{\mathrm{H}_1}}
\newcommand{\hn}{\ensuremath{\mathrm{H}_n}}

\newcommand{\z}{\ensuremath{\mathbb{Z}}}

\newcommand{\hnj}{\ensuremath{\mathrm{H}_{nj}}}

\newcommand{\gzi}{\ensuremath{\mathrm{G}_{0i}}}
\newcommand{\goi}{\ensuremath{\mathrm{G}_{1i}}}
\newcommand{\gni}{\ensuremath{\mathrm{G}_{ni}}}

\newcommand{\hzl}{\ensuremath{\mathrm{H}_{0j}}}
\newcommand{\hol}{\ensuremath{\mathrm{H}_{1j}}}

\newcommand{\ZZ}{\mathbb{Z}}
\newcommand{\NN}{\mathbb{N}}
\newcommand{\QQ}{\mathbb{Q}}
\newcommand{\CC}{\mathbb{C}}
\newcommand{\im}{\operatorname{im}}
\newcommand{\myTor}[2][n]{{#2}[#1]}

\newcommand{\cstar}{\mbox{$C^*$}}
\newcommand{\K}{\mathrm{K}}
\newcommand{\pK}[2][n]{\mathrm{K}_0({#2};\ZZ/{#1})}
\newcommand{\nzcc}[1][n]{$#1$-coefficient complex\xspace}
\newcommand{\nzccs}[1][n]{$#1$-coefficient complexes\xspace}
\newcommand{\CS}{C(\mathrm{S}^1)}
\newenvironment{demo}{{\bf\sl Proof:}}{\hfill$\Box$}

\begin{document}

\title[Inductive Limits of $\mathrm{K}$-theoretic Complexes]{Inductive Limits of $\mathrm{K}$-theoretic Complexes with
Torsion Coefficients}
\author{S\o ren Eilers} 
\address{Department of Mathematics, University of Copenhagen, Universitetsparken 5, 2100 K\o benhavn, Danmark}
\email{eilers@math.ku.dk}
\author{Andrew S. Toms}
\address{Department of Mathematics and Statistics, University of New Brunswick, Fredericton, New Brunswick, Canada, E3B 5A3}
\email{atoms@unb.ca}

\begin{abstract}
We present the first range result for the total $\mathrm{K}$-theory
of $C^*$-algebras.  This invariant has been used successfully to
classify certain separable, nuclear $C^*$-algebras of real rank zero.
Our results complete the classification of the so-called AD algebras
of real rank zero.
\end{abstract}

\maketitle

\section{Introduction}

A theorem which classifies the objects of a category up to some notion
of equivalence via an invariant begs naturally the question of range for
the classifying invariant.  In the
classification theory of \cstar-algebras by $\K$-theoretic invariants the fundamental
range-of-invariant result is the theorem of Effros, Handelman, and Shen
(\cite{ehs}), which states that the ordered groups arising as
$\kz$-groups of AF algebras are exactly the dimension groups 
studied first by Riesz (\cite{fr:qnftgol}) and Fuchs.
(\cite{lf:rg}).  Post AF classification results for $C^*$-algebras have
required invariants more complex than $\kz$ alone, yet it has typically been
possible to pair such results with an Effros-Handelman-Shen type theorem 
establishing the range of the classifying invariant.  Notable examples include \cite{jv:ah}
and \cite{gae:dgt}. 

The aim of the present paper is to give the first range-of-invariant result
associated to the classification of certain \cstar-algebras
of real rank zero, completed by Dadarlat and Gong in \cite{mdgg:crahcrrz}.
To obtain complete invariants for this  class of non-simple
\cstar-algebras, the ordered $\K_*$-group $\kz(-)\oplus\ko(-)$ must be
augmented:  the addition of ordered $\K$-groups with torsion coefficients and 
certain natural homomorphisms between them is required.  While the completeness of this invariant 
has been established for almost a decade in a number of cases 
(\cite{se:ciabtk}, \cite{mdtal:umctkg}, \cite{mdgg:crahcrrz}),
there have been no {range results} available until now.

The situation is complicated by the intricate nature of the
augmented $\K$-theory. Following \cite{mdgg:crahcrrz}, this invariant
associates to each \cstar-algebra $\A$ a family  of groups
\[
\K_0(\A), \K_1(\A), \pK{\A}, \pKone{\A}
\]
with $n$ ranging over $\{2,3,\dots\}$, as well as an order structure on
\[
\K_0(\A)\oplus \K_1(\A)\oplus\bigoplus_{n\in\{2,3,4,\dots\}}
\left[\pK{\A}\oplus \pKone{\A}\right]
\]
and families of group homomorphisms
\begin{gather*}
\rho^i_n:\K_i(\A)\longrightarrow \pKi{\A}\\
\beta^i_n:\pKi{\A}\longrightarrow \K_{i+1}(\A)\\
\kappa^i_{n,m}:\pKi{\A}\longrightarrow\pKi[m]{\A}.
\end{gather*}
Thus, an isomorphism of invariants amounts to a family
\[
\phi_i:\K_i(\A)\longrightarrow \K_i(\B)\qquad
\psi_i:\pKi{\A}\longrightarrow \pKi{\B}
\]
which preserves the order structure and intertwines all morphisms
$\rho,\beta,\kappa$.

 To keep technicalities to a minimum while
staying in a class where, by the counterexamples given in 
\cite{mdtal:ccomk} and \cite{mdse:bmn}, the full force of such an invariant is
really needed, we
shall concentrate on the class of so-called AD algebras of real rank
zero.  Recall that an AD algebra is an inductive limit of finite direct sums
of matrix algebras over elements of
\[
\mathrm{D} :=\{\CC,C(S^1),\mathbb{I}_2^{\sim},\mathbb{I}_3^{\sim},\mathbb{I}_4^{\sim},\dots\}, 
\]
where
$\mathbb{I}^ {\sim}_n$ is the dimension drop algebra
\[
\{f\in C([0,1],\mathbf{M}_n(\CC))\mid f(0),f(1)\in\CC\mathbf{1}\}.
\]
Such \cstar-algebras may be classified by a more
manageable invariant of the form
\[
\K_0(\A)\otimes\QQ\longrightarrow \K_0(\A;\QQ/\ZZ)\longrightarrow \K_1(\A)
\]
provided that they are
of real rank zero, cf. \cite{mdse:ccpikc}. Here, as we shall recall
below, $\K_0(\A;\QQ/\ZZ)$ should be thought of as a kind of
conglomerate of $\pK{\A}$ for all $n$. If the torsion part
of
$\K_1(\A)$ is annihilated by a fixed integer $n$, then the even more
manageable
invariant
\[
\BK(\A)\qquad \K_0(\A)\longrightarrow \pK{\A}\longrightarrow \K_1(\A)
\]
suffices.

Our strategy is to establish a range result in the latter case
first, then use this to derive the general result.  The key technical
element of the proof is a decomposition result for refinement monoids
attributed to Tarski by Wehrung (\cite{fw:ipomi}).

To illustrate our results we revisit examples of AD algebras
originally considered by Dadarlat and Loring, which showed
that such algebras could have isomorphic ordered $\kx$-groups without having
isomorphic augmented K-theory. Using our main result, we
parametrize the real rank zero AD algebras with $\kx$-groups as
considered in \cite{mdtal:ccomk}, and show that there are uncountably
many non-isomorphic such algebras.

\section{Building blocks}

In this section we introduce the notion of an \emph{n-coefficient
complex}.  This type of object is meant to abstract the characteristics
of certain augmented K-theoretic invariants for AD algebras of real 
rank zero --- invariants which will be reviewed in detail in section 3.
We begin with some preliminaries and notation.

A \emph{graded ordered group} is a graded group $G_0\oplus G_1$ in
which the $G_0$-component dominates the order in the sense that
\[
\left.\begin{array}{c}(x,y)\geq 0\\(x,y')\geq 0\end{array}\right\}
\Longrightarrow (x,y\pm y')\geq 0
\]
For any group $G$, we denote by $G[n]$ the subgroup of elements of $G$
annihilated by $n\in\NN$. When $G$ is an ordered group, we
denote by $\Ideal(x)$ the order ideal containing $x$.
Recall the notions of unperforated and weakly unperforated groups from
\cite{krg:poagi}.

Let $G$ be an ordered abelian group, $H$ an abelian group, and $f:G \rightarrow H$ 
a surjective group homomorphism.  Say that $h \in H$ is positive if it is the image
of a positive element in $G$.  The important and obvious feature of the order on $H$ thus 
defined (the so-called \emph{quotient order}) is that every positive element in $H$
lifts to a positive element in $G$ (\cite{krg:poagi}).

\noindent
\begin{defs}\label{nccdef} Let $n\in\{2,3,\dots\}$. An $n$-coefficient complex $\overline{G}$ is an exact sequence
\[
\gz \stackrel{\rho}{\longrightarrow} \gn
\stackrel{\beta}{\longrightarrow} \go
\] 
of abelian groups which, setting
\[
\gx:= \gz \oplus \go, \ \ \gd := \gz \oplus \gn,
\]
has the following properties:
\begin{enumerate}[(i)]
\item $n\gn=0$
\item $\ker\rho=n \gz, \im\beta=\myTor{\go}$.
\item $\gx$ and $\gd$ are graded ordered
groups restricting to the same order on $\gz$
\item $\gx$ has the Riesz interpolation property.
\item $\gz \oplus \rho(\gz)$ has the quotient order coming from $\mathrm{id}_{G_0} \oplus \rho$
\item $\gz \oplus \beta(\gn)$ has the quotient order coming from $\mathrm{id}_{G_0} \oplus \beta$ 
\item $\gz$ is unperforated and $\gx$ is weakly unperforated.
\end{enumerate}
We say that an element $(x,y,z)$ is positive in
$\overline{G}$ if and only if 
\[
gd \ni (x,y) \geq 0 and \gx \ni (x,z) \geq 0.
\] 
\end{defs}
\vspace{3mm}

A morphism $\overline{\theta}: \overline{G} \rightarrow \overline{H}$ of 
\nzccs is a positive ordered triple of
 linear maps $(\theta_0, \theta_n, \theta_1)$ such that
\[
\theta_0:\gz \rightarrow \hz,\ \theta_n:\gn \rightarrow \hn,\ 
\theta_1:\go \rightarrow \ho,
\]
and the maps commute with $\rho$ or $\beta$ as appropriate.

We conclude this section by introducing three types of $n$-coefficient complexes --- our so-called
building blocks.

\vspace{3mm}
\noindent
\begin{defs}\label{bbs} 
The
\begin{enumerate}[(i)]
\item $(\mathbb{C},n)$ complex
\begin{gather*}
\z \stackrel{\rho}{\longrightarrow} \z/n
\stackrel{\beta}{\longrightarrow} 0\\
\rho:1\mapsto \overline{1}, \beta:\overline{1}\mapsto 0,
\end{gather*}
\item $(\mathbb{I}^{\sim}_m,n)$ complex 
\begin{gather*}
\z \stackrel{\rho}{\longrightarrow} \z/n\oplus \z/m
\stackrel{\beta}{\longrightarrow} \z/m\\
\rho:1\mapsto (\overline{1},\overline{0});
\beta:(\overline{1},\overline{0})\mapsto 0,
(\overline{0},\overline{1})\mapsto \overline{1}
\end{gather*}
and
\item $(\CS,n)$ complex 
\begin{gather*}
\z \stackrel{\rho}{\longrightarrow} \z/n
\stackrel{\beta}{\longrightarrow} \z\\
\rho:1\mapsto \overline{1}, \beta:\overline{1}\mapsto 0,\\
\end{gather*}
\end{enumerate}
where \gd\ and \gx\ have the strict order coming from the first direct summand,
are $n$-coefficient complexes.
\end{defs}

The motivation for these defining these objects, hinted at by 
their very names, will be made clear in the following section.


\section{$\mathrm{K}$-theory with Coefficients}

In this section we collect a suite of known results which
together prove that \nzccs appear as the $\K$-theory of certain
\cstar-algebras. 

The following definitions, originating in the work of Dadarlat, 
Gong, and the first named author (see \cite{mdse:ccpikc})
 are based on the observation (from
\cite[2.3]{sz:rdpisma}) that the lattices of order ideals of
$\K_0(\A)$ and of ideals of $\A$ are naturally isomorphic for
\cstar-algebras with minimal ranks.

\begin{defs}\label{egsemap}
Let $\A$ be a \cstar-algebra of real rank zero and stable rank
one. When $I$ is an order ideal of $\K_0(\A)$ we define 
\[
\egmap{I}=\iota_*(\K_1({\mathfrak{I}})) \ \ \mathrm{and} \ \ \semap{I}=\iota_*(\pK{\mathfrak{I}}),
\]
where $\mathfrak{I}$ is the unique ideal of $\A$ with
$I=\iota_*(\K_0(\mathfrak{I}))$ and $\iota:\mathfrak{I}\hookrightarrow \A$ is the inclusion map.
\end{defs}

We now equip $\kx(\A)=\K_0(\A)\oplus \K_1(\A)$ and
$\kd(\A)=\K_0(\A)\oplus\pK{\A}$ with the orders given by 
\[
(x,y)\geq 0\Longrightarrow\left\{
\begin{array}{c}
x\geq 0\\
y\in\egmap{\Ideal(x)}\end{array}\right\}
\]
and
\[
(x,z)\geq 0\Longrightarrow\left\{
\begin{array}{c}
x\geq 0\\
z\in\semap{\Ideal(x)}\end{array}\right\},
\]
respectively.

It is well known (cf.\ \cite{segae:rpkacmr}) that the order thus defined on $\kx(\A)$ will
coincide with the standard order on $\kx(\A)$ derived from the
isomorphism
\[
\kx(\A)\cong \kz(\A\otimes C(S^1))
\]
In general, the
order on $\kd(\A)$ will \underline{not} be the one similarly derived
from the isomorphism
\[
\kd(\A)\cong \kz(\A\otimes \mathbb{I}^ {\sim}_n).
\]
But since, as seen in
\cite{mdse:ccpikc}, these
two order structures
allow the same positive group isomorphisms for a large class of
\cstar-algebras including the $AD$ algebras, the choice of order structure for 
the invariant has no influence on the associated classification results.

\begin{props}\label{isnzcc}
Let $\A$ be a \cstar-algebra of real rank zero and stable rank one.
Assume that $\K_*(\A)$ is weakly unperforated, and that $\K_0(\A)$ is unperforated. For any $n \in \{2,3,\ldots\}$,
\[\xymatrix{{\BK(\A):}&
{\K_0(\A)}\ar[r]&
{\K_0(\A;\ZZ/n)}\ar[r]&
{\K_1(\A)}
}
\]
is an $n$-coefficient complex.
\end{props}
\begin{demo}
We verify properties $(i)-(vii)$ from Definition \ref{nccdef}.

The purely algebraic properties $(i)$ and $(ii)$ hold true for any such
sequence, cf.\ \cite{cs:tmcIV}. Furthermore, it is clear from our
definition of the order on $\K_*(\A)$ and $\kd(\A)$ that they are graded order groups
based on the same order on $\K_0(\A)$. Since we have noted that we are
in fact working with the standard order on $\K_*(\A)$, \cite{lgb:ripk}
or \cite{segae:rpkacmr} 
show that condition $(iv)$ is met.

Inspection of the diagram
\begin{gather}\label{tochase}
\xymatrix{
{\K_0(\mathfrak{I})}\ar[r]^-{\rho}\ar[d]_-{\iota_*}&
{\pK{\mathfrak{I}}}\ar[r]^-{\beta}\ar[d]_-{\iota_*}&
{\K_1(\mathfrak{I})}\ar[d]_-{\iota_*}\\
{\K_0(\A)}\ar[r]_-{\rho}&
{\pK{\A}}\ar[r]_-{\beta}&
{\K_1(\A)}
}
\end{gather}
when $\mathfrak{I}$ is the ideal of $\A$ corresponding to the order ideal $\Ideal(x)$ shows that
if $\gz \ni x \geq 0$, then $\gn \ni (x,\rho(x)) \geq 0$, and, similarly, that if 
$\gn \ni (x,y) \geq 0$, then $\gx \ni (x,\beta(y)) \geq 0$.

To prove
property $(v)$ we look again at the diagram \eqref{tochase}. By
assumption, $y\in\pK{\A}$ is in the image of both maps with target
$\pK{\A}$, so an easy diagram chase gives the desired result whenever 
$\iota_*:\K_1(\mathfrak{I})\longrightarrow \K_1(\A)$ is injective. But $\iota_*$ is
always injective by \cite{hlmr:eilca}.  Combining this fact with the observation 
of the preceding paragraph, we have property $(v)$.

For $(vi)$, we do a similar diagram chase.
Finally, we have explicitly required the properties in $(vii)$.
\end{demo}



\section{Decomposition Lemmas}

In this section we establish some decomposition results in the spirit of Riesz 
for $n$-coefficient complexes.  These lemmas will allow us to prove 
an Effros-Handelman-Shen-type result for these complexes, realising them as inductive
limits of our building blocks.

We shall rely heavily on results in \cite{gae:dgt} pertaining to the
family of ordered $\kx(-)$-groups of AH algebras
with real rank zero. These groups have the following property:

\begin{defs}\label{wup} (Cf. Goodearl (\cite[Lemma 8.1]{krg:kmacrrz}) and
Elliott (\cite{gae:dgt})) An ordered group $G$ is said to be
\emph{weakly unperforated} if
\begin{enumerate}[(i)]
\item whenever $mx\in G_+$ there exists $t\in \tor(G)$ with $x+t\in G_+$ and
$mt=0$;
\item whenever $y\in G_+$, $t\in\tor(G)$, and $ny+t\in G_+$ for some
  $n\in\NN$, then $y\pm t\in G_+$.
\end{enumerate}
\end{defs}
Note that property $(ii)$ is automatic in our case since all torsion is
localized in the odd part of a graded ordered group.  Although we do not apply
the next observation in the sequel, we nevertheless record it for possible
future use:  all of the results in
this section hold true if the condition of unperforation in $\gz$
in Definition \ref{nccdef}$(vii)$ is relaxed to weak unperforation. 

\begin{lms}[{Elliott (\cite[Corollary 6.6]{gae:dgt})}]\label{fourone}
Let $\gx = \gz \oplus \go$ be a weakly unperforated graded ordered group
with the Riesz decomposition property.  If $s_1,\ldots,s_m \leq g$ 
where $g \in \gz^+$ and $s_1,\ldots,s_m \in \go$, then $g = g_1 +
\cdots + g_m$ with $g_1,\ldots,g_m \in \gz^+$ and $s_i \leq g_i$.
\end{lms}

We say that a family $H_1,\dots,H_n$ of subgroups of a given groups
$G$ is \emph{independent} if
\[
\sum_{i=1}^nx_i=0,x_i\in H_i\Longrightarrow x_1=\cdots=x_n=0
\]

\begin{lms}[{Elliott (\cite[Corollary 6.3]{gae:dgt})}]\label{fourtwo}
Let $\gx = \gz \oplus \go$ be a weakly unperforated graded ordered
group
with the Riesz decomposition property.  Suppose $x \leq \sum_{j=1}^{k}
g_j$, where $x \in \go$ and $g_j \in \gz^+$.  Then,
there exist an independent family $H_j$, $j \in \{1,\ldots,k\}$, of 
finitely generated subgroups of \go\ such that $H_j \leq g_j$ and a decomposition
$x=\sum_{j=1}^{k} x_j$ such that $x_j \in H_j$.
\end{lms}

Note that if $x$ is as in Lemma \ref{fourtwo} and has order $m$, then each $x_j$
has order at most $m$ by the independence of the $x_j$. Thus, by property $(ii)$ of
Definition \ref{nccdef}, if the \gx\ of Lemma \ref{fourtwo} is in fact \gx\ for
some $n$-coefficient complex $\overline{G}$ and $x$ is in the image of $\beta$,
then so too are the $x_j$.

Wehrung attributes the following observation to Tarski:

\begin{lms}[{Cf.\ Wehrung (\cite[Lemma 1.9]{fw:ipomi})}]\label{fourthree}
Let $\gz$ be an ordered group with the Riesz interpolation property, and let
$a,b\in\gz$ satisfy $a,b\geq 0$ and $a \leq nb$, $n \in \NN$.  Then, there
exist $b_0,\ldots,b_n\geq 0$ such that $b = \sum_{i=0}^n b_i$ and $a =
\sum_{i=1}^n i b_i$.
\end{lms}

\begin{lms}\label{fourfour}
Let $n\in\{2,3,\dots\}$ and let $\overline{G}$ be an $n$-coefficient complex.
Let $(e,f,g) \in \overline{G}$ be a positive element and let 
there be given a decomposition $e = \sum_{j=1}^k e_j$, $e_j \in
\gz^+$.  Then, there exist elements
\[
g_1,g_2,\ldots,g_k \in \go \ \ \mathrm{and} \ \  f_1,f_2,\ldots,f_k \in \gn
\]
such that 
\[
(e,f,g) = \sum_{j=1}^{k} (e_j,f_j,g_j)
\]
and $(e_j,f_j,g_j)$ is positive in $\overline{G}$ for each $j \in
\{1,\ldots,k\}$.
\end{lms}

\begin{demo}
Since $\beta(f), g \leq \sum_{j=1}^{k} e_j$, there exist elements
\[
g_1,\ldots,g_k,l_1,\ldots,l_k \in \go
\]
such that
$\beta(f)=\sum_{j=1}^k l_j$ and $g=\sum_{j=1}^k g_j$ with $l_j,g_j
\leq e_j$ (Lemma \ref{fourtwo}). As noted in the comment following
that lemma, we may assume that $l_i\in\im\beta$, so by $(vi)$ of Definition 2.1, 
the $l_j$ have $\beta$-lifts $\widetilde{l_j}$ such that
$\widetilde{l_j} \leq e_j$.  Thus, both $f$ and $\sum_{j=1}^k \widetilde{l_j}$
are majorised by $e$ and have the same image under $\beta$.  We
conclude that the difference $f - \sum_{j=1}^k \widetilde{l_j}$ is in the image of
$\rho$, and is majorised by $e$.  By property $(v)$ of Definition 2.1, we may choose
$c\in \Ideal(e)$  so that $f = \sum_{j=1}^k
\widetilde{l_j} + \rho(c)$, $c \in \mathrm{I}(e)$.  Since $\mathrm{I}(e) =
\mathrm{I}(e_1) + \cdots + \mathrm{I}(e_k)$, there is a decomposition
$c = c_1+\cdots+c_k$, $c_j \in \mathrm{I}(e_j)$.  Put $f_j =
\widetilde{l_j} + \rho(c_j)$, so that $f= \sum_{j=1}^k f_j$.  By
construction, we have $f_j,g_j \leq
e_j$, so that $(e_j,f_j,g_j)$ is positive in the $n$-coefficient complex
for each $j \in \{1,\ldots,k\}$.  \end{demo} 

\vspace{3mm}
\noindent
Note that the lemma above holds even when one specifies the $g_j \leq
e_j$ \textit{a priori}.

In the following, we will use the term \emph{refinement} of a
collection of elements $x_1,\dots,x_s$ to denote a new collection of
elements $\widetilde{x}_1,\dots,\widetilde{x}_t$ with the property that
$\{1,\dots,t\}$ can be partitioned into $s$ subsets, such that the sum
of the elements corresponding to the indices in the $j$th subset is
exactly $x_j$.

\begin{lms}\label{fourfive}
Fix $n\in\{2,3,\dots\}$ and let $\overline{G}$ be an  $n$-coefficient complex.  Let 
$(e_i,f_i,0)$, $i \in \{1,\ldots,k\}$, be positive in $\overline{G}$. 
Let there be given elements $x_1,\ldots,x_r \in \gz^+$
and $z_1,\ldots,z_r \in \go$ such that $z_j \leq x_j$, 
and non-negative integers $\lambda_{ij},
\delta_{ij}$, $i \in \{1,\ldots,k\}$, $j \in \{1,\ldots,r\}$, such
that
\[
e_i = \sum_{j=1}^r \lambda_{ij}x_j, \ \ \beta(f_i) = \sum_{j=1}^r
\delta_{ij} z_j.
\]
Then, there exist refinements $\widetilde{x_1},\ldots,\widetilde{x_s}$ of $x_1,\ldots,x_r$
and $\widetilde{z_1},\ldots,\widetilde{z_s}$ of $z_1,\ldots,z_r$, and lifts
$\widetilde{y_l} \in \gn$ of the $\widetilde{z_l}$ with $\widetilde{y_l} \leq
\widetilde{x_l}$, having the following property:
there are non-negative integers $\gamma_{il}$, $\kappa_{il}$,
and $n_{il}$, $i \in \{1,\ldots,k\}$ and $l
\in \{1,\ldots,s\}$, such that
\[
e_i = \sum_{l=1}^s \gamma_{il} \widetilde{x_l}
\]
and
\[
f_i = \sum_{l=1}^s \kappa_{il} \widetilde{y_l} + n_{il} \rho (\widetilde{x_l}).
\]
Furthermore, $\gamma_{il} \neq 0$ whenever $n_{il} \neq 0$.
\end{lms}

\begin{demo} Following the proof of Lemma \ref{fourfour}, we may assume that we have lifts
$y_j$ of each $z_j$, and positive elements $c_{ij} \in
\mathrm{I}(x_j)$  such that 
$f_i = \sum_{j=1}^r \delta_{ij} y_l +\rho(c_{ij})$ and $\delta_{ij}
y_j + \rho(c_{ij}) \leq  x_j$, $1 \leq j \leq r$, $1 \leq
i \leq k$.  Fix $j$.  By Lemma \ref{fourthree} there is, for each $i \in
\{1,\ldots,k\}$, a decomposition $x_j = x^i_{j,1}+\cdots+x^i_{j,k_i}$,
$k_i \in \NN$,
such that $c_{ij}$ is in the non-negative integral linear span of
$\{x^i_{j,1},\ldots,x^i_{j, k_i}\}$.  Choose by the Riesz
 property in $\gz$ a decomposition $x_j =
x_{j,1}+\cdots+ x_{j, m_j}$, some $m_j \in \NN$, which
simultaneously refines all of the $x_j = x^i_{j,1}+\cdots+x^i_{j, k_i}$
decompositions, $1 \leq i \leq k$.  Then, there exist non-negative
integers $n^i_{j,1},\ldots,n^i_{j, m_j}$, $1 \leq i \leq k$, such that
$c_{ij} = n^i_{j,1} x_{j,1} + \cdots + n^i_{j, m_j} x_{j, m_j}$.
Since $\delta_{ij}y_j + \rho(c_{ij}) \leq  x_j$ we have that $\delta_{ij}
y_j  \leq  x_j$ (property $(v)$ of Definition 2.1), whence $y_j \leq x_j = x_{j,1}
+ \cdots + x_{j, m_j}$ (Lemma \ref{fourone}).  By
Lemma \ref{fourfour} there is a decomposition $y_j = y_{j,1} + \cdots + y_{j, m_j}$
such that $y_{j,l} \leq x_{j,l}$, $1 \leq l \leq m_j$.  Thus, 
\[
\delta_{ij} y_j + \rho(c_{ij}) = \sum_{p=1}^{m_j} \delta_{ij} y_{j,p}
+ n^i_{j,p} \rho(x_{j,p}).
\]
Define 
\[
\{\widetilde{x_1},\ldots,\widetilde{x_s}\} := \cup_{j \leq r}
\{x_{j,1},\ldots,x_{j, m_j}\},
\]
\[
\{\widetilde{y_1},\ldots,\widetilde{y_s}\} := \cup_{j \leq r}
\{y_{j,1},\ldots,y_{j, m_j}\},
\]
and $\widetilde{z_l} := \beta(\widetilde{y_l})$.  The lemma follows. \end{demo}

\section{Bounded Torsion}
\subsection{A Local Property}
In this section we establish that 
\nzccs satisfy a local property such as the one
whose importance was realized by Shen (cf.\ \cite{ehs}) 
in the setting of classical dimension groups.

\begin{lms}\label{fiveone}
Let $\overline{\mathcal{G}}$ be an $n$-coefficient complex. 
Let $\overline{G} = \oplus_{i=1}^n \overline{G_i}$ be a
direct sum of $n$-coefficient complex (with the direct sum order
structure), where each $\overline{G_i}$ is a
$(\mathbb{C},n)$, $(\mathbb{I}^{\sim}_m,n)$ or $(\CS,n)$
complex.  Let $\overline{\theta}: \overline{G} \rightarrow \overline{\mathcal{G}}$ be
a morphism.  Then, there exist an $n$-coefficient complex $\overline{H} =
\oplus_{j=1}^m \overline{H_j}$ with each $\overline{H_j}$ a $(\mathbb{C},n)$, 
$(\mathbb{I}^{\sim}_m,n)$ or $(\CS,n)$ complex, and morphisms
$\overline{\gamma}: \overline{G} \rightarrow \overline{H}$ and $\overline{\lambda}:
\overline{H} \rightarrow \overline{\mathcal{G}}$ such that the diagram
\[
\xymatrix{
{\overline{G}}\ar[dr]^-{\overline{\gamma}}\ar[dd]^-{\overline{\theta}}&\\
&{\overline{H}}\ar[dl]^-{\overline{\lambda}}\\
{\overline{\mathcal{G}}}&}
\] 
commutes and $\ker\overline\gamma=\ker\overline\theta$.
\end{lms}

\begin{demo} Suppose that the conclusion above is relaxed to read 
\[
\ker \overline{\gamma} = \ker \overline{\theta} \ 
(\mathrm{mod} \rho(\gz)),
\]
with all other things being equal.  Then, the original conclusion of the
lemma follows.  Indeed, suppose that $a \in \ker
\overline{\theta}$ is in the image of $\rho$, i.e., $a = \rho(q)$ for some
$q \in \gz^+$.  Then, $\theta_0(q) =  n \cdot q^{'}$
for some $q^{'} \in \mathcal{G}_0^+$.  Put $\overline{H}^{'} = \overline{H} \oplus
\overline{R}$, where $\overline{R}$ is a $(\mathbb{C},n)$ complex, and extend
$\overline{\lambda}$ to $\overline{H}^{'}$ by sending the positive generator of 
$R_0$ to $q^{'}$.  Apply the weakened conclusion above to find a
direct sum of building block complexes $\overline{H}^{''}$ and maps 
$\overline{\gamma}^{''}: \overline{H}^{'} \rightarrow \overline{H}^{''}$ and  $\overline{\lambda}^{''}:
\overline{H}^{''} \rightarrow \overline{\mathcal{G}}$ such that the diagram
\[
\xymatrix{
{\overline{G}}\ar[dr]^-{\overline{\gamma}}\ar[ddd]^-{\overline{\theta}}&&\\
&{\overline{G}'}\ar[dr]^-{\overline{\gamma}''}\ar[ldd]^-{\overline{\theta}'}&\\
&&{\overline{H}''}\ar[lld]^-{\overline{\lambda}''}&\\
{\overline{\mathcal{G}}}&&}
\] 
commutes.  Note that $\ker( \overline{\gamma}^{''} \circ \overline{\gamma} )
= \ker \overline{\theta} \ (\mathrm{mod} \rho(\gz))$.  Furthermore,
$a \in \ker (\overline{\gamma}^{''} \circ \overline{\gamma})$, since
$(\gamma_0^{''} \circ \gamma_0)(q)$ must be $\mathrm{ord}(\rho(q))$ times 
the image of the positive generator of $R_0$ under $\gamma_0^{''}$.
Repeating this procedure for each of the finitely many elements
in $\ker \overline{\theta} \cap \rho(\gz)$ yields the conclusion of
the lemma proper.

It remains to prove that the lemma holds if we only require that
$\ker \overline{\gamma} = \ker\overline{\theta} \ 
(\mathrm{mod} \rho(\gz))$.  Let $\gzi = \langle e_i
\rangle$ ($e_i \geq 0$), $\goi = \langle g_i \rangle$, and choose $f_i \in \gni$
(necessarily $\leq e_i$) such that $\gni =\langle \rho(e_i) \rangle
\oplus \langle f_i \rangle$.  Note that if $\overline{G}_i$ is a
$(\mathbb{C},n)$ or $(\CS,n)$ complex, then we may (and do) take $f_i = 0$.  In
the case of a $(\mathrm{I}_m,n)$ complex, $\beta(f_i) = g_i$.  Define
$a_i:= \theta_0(e_i)$, $b_i:= \theta_n(f_i)$ and $c_i:=
\theta_1(g_i)$.  By the main Theorem in Section 5 of \cite{gae:dgt}, there is a
complex $\overline{H}^{'} = \oplus_{l=1}^k \overline{H}^{'}_l$ with each
$\overline{H}^{'}_l$ a building block (and elements $e^{'}_l$, $f^{'}_l$
and $g^{'}_l$ playing roles analogous to those of the $e_i$, $f_i$ and
$g_i$ above), and maps
\[
\gamma_0^{'}:\gz \rightarrow \hz^{'},\ \gamma_1^{'}:\go \rightarrow
\ho^{'}
\]
and
\[
\lambda_0^{'}:\hz^{'} \rightarrow \mathcal{G}_0,\ \lambda_1^{'}:
\ho^{'} \rightarrow \mathcal{G}_1
\]
such that 
\[
\xymatrix{
{\gz}\ar[rd]^-{\gamma_0^{'}}\ar[dd]\ar[rr]^-{\rho}&&{\gn}\ar'[d][dd]\ar[rr]^-{\beta}&&{\go}\ar'[d][dd]\ar[rd]^-{\gamma_1^{'}}&\\
&{\hz^{'}}\ar[ld]^-{\lambda_0^{'}}\ar[rr]^(.3){\rho}&&{\hn^{'}}\ar[rr]^(.3){\beta}&&{\ho^{'}}\ar[ld]^-{\lambda_1^{'}}\\
{\mathcal{G}_0}\ar[rr]^-{\rho}&&{\mathcal{G}_n}\ar[rr]^-{\beta}&&{\mathcal{G}_1}&}
\]
commutes and $\ker(\theta_0, \theta_1) =
\ker(\gamma_0^{'},\gamma_1^{'})$.  To be fair, \cite{gae:dgt} only
provides the $\hzl^{'}$ and $\hol^{'}$, but we may clearly associate a
building block complex to any pair $(\hzl^{'},\hol^{'}) = (\ZZ,R)$, $R \in \{0,
\ZZ, \ZZ/2, \ZZ/3,\dots \}$.  
Of course, this association is token for
now, as we have no positive maps $\gamma_n^{'}:\gn \rightarrow
\hn^{'}$ and $\lambda_n^{'}:\hn^{'} \rightarrow \mathcal{G}_n$ to fill
in the diagram above.

We have $a_i = \sum_{l=1}^k \kappa_{il} \lambda_0^{'}(e_l^{'})$,
$\ZZ \ni \kappa_{il} \geq 0$, and $\beta(b_i) = \sum_{l=1}^k \phi_{il}
\lambda_1^{'}(g_l^{'})$, $\phi_{il} \in \ZZ$.  By Lemma \ref{fourfive}
there exist refinements 
\[
\{\widetilde{a_1},\ldots,\widetilde{a_m}\} \ (\widetilde{a_j}
\geq 0, 1 \leq j \leq m) \ \mathrm{and} \
\{\widetilde{c_1},\ldots,\widetilde{c_m}\}
\]
of 
\[
\{\lambda_0^{'}(e_1^{'}),\ldots,\lambda_0^{'}(e_k^{'})\} \ \mathrm{and} \
\{\lambda_1^{'}(g_1^{'}), \ldots, \lambda_1^{'}(g_k^{'})\},
\]
respectively, and elements
$\widetilde{b_1},\ldots,\widetilde{b_s} \in \mathcal{G}_n$ such that for some integers
$\zeta_{ij}$, $\xi_{ij}$, $\iota_{ij}$ and $n_{ij}$, $1 \leq i \leq k$, $1 \leq j
\leq m$, we have
\[ 
a_i = \sum_{j=1}^m \zeta_{ij} \widetilde{a_j},
\]
\[
b_i = \sum_{j=1}^m \xi_{ij} \widetilde{b_j} + n_{ij} \rho (\widetilde{a_j}),
\]
and
\[ 
c_i = \sum_{j=1}^m \zeta_{ij} \widetilde{c_j}.
\]
Furthermore, $\widetilde{b_j}$ is a lift of $\widetilde{c_j}$ whenever
$\widetilde{c_j}$ is in the image of $\beta$, and is zero otherwise.  Note
that the $\widetilde{c_j}$ can be and should be chosen to be in the 
image of $\beta$ whenever it is a torsion element.
By the main Theorem of Section 5, \cite{gae:dgt}, there is an $n$-coefficient complex
$\overline{H}= \oplus_{j=1}^m \overline{H}_j$, $\overline{H}_j$ a building block
complex for each $j$, and there are maps
\[
\gamma_0:\gz \rightarrow \hz,\ \gamma_1:\go \rightarrow \ho
\]
and
\[
\lambda_0:\hz \rightarrow \mathcal{G}_0,\ \lambda_1:
\ho \rightarrow \mathcal{G}_1
\]
such that (with the dotted arrows representing desired but as yet
undefined maps)
\[
\xymatrix{
{\gz}\ar[rd]^-{\gamma_0}\ar[dd]\ar[rr]^-{\rho}&&{\gn}\ar'[d][dd]\ar[rr]^-{\beta}\ar@{-->}[rd]&&{\go}\ar'[d][dd]\ar[rd]^-{\gamma_1}&\\
&{\hz}\ar[ld]^-{\lambda_0}\ar[rr]^(.3){\rho}&&{\hn}\ar[rr]^(.3){\beta}\ar@{-->}[ld]&&{\ho}\ar[ld]^-{\lambda_1}\\
{\mathcal{G}_0}\ar[rr]^-{\rho}&&{\mathcal{G}_n}\ar[rr]^-{\beta}&&{\mathcal{G}_1}&}
\]
commutes, and $\gamma_0$ and $\gamma_1$ factor
through $\hz^{'}$ and $\ho^{'}$, respectively.  Thus, $\ker(\theta_0, \theta_1) =
\ker(\gamma_0,\gamma_1)$.  Let $\widetilde{e_j}$,
$\widetilde{f_j}$, and $\widetilde{g_j}$ play roles in $\overline{H}_j$ analogous
to the roles of the $e_i$, $f_i$, and $g_i$ in $\overline{G}_i$.  Then 
$\lambda_0(\widetilde{e_j})=\widetilde{a_j}$ and
$\lambda_1(\widetilde{g_j})=\widetilde{c_j}$.  
For every pair $(i,j)$, $1 \leq i \leq k$, $1\leq j \leq m$, 
define a partial map $\gamma_n^{ij}:\gni \rightarrow \hnj$
by 
\[      
\gamma_n^{ij}(f_i) := \xi_{ij} \widetilde{f_j} + n_{ij} \rho(\widetilde{e_j}).
\]  
Let $\lambda_n:\hn \rightarrow \mathcal{G}_n$ be defined by
$\lambda(\widetilde{f_j}) := \widetilde{b_j}$, and put $\gamma_n = \oplus_{j=1}^m
(\sum_{i=1}^k \gamma_n^{ij})$.  These maps are positive and so complete
the morphisms $\overline{\lambda}$ and $\overline{\gamma}$, establishing the desired weak
version of the lemma. \end{demo}

\subsection{The Range Result}

\begin{lms}\label{hit}
Fix $n\in\{2,3,\dots\}$,  let $\overline{\mathcal{G}}$ be a $n$-coefficient complex, 
and consider a positive element
$(e,f,g)\in\overline{\mathcal{G}}$. There exists a $n$-coefficient complex
$\overline{H}$ which is a finite direct sum of $(\mathbb{C},n)$,
$(\mathbb{I}^{\sim}_m,n)$ and $(\CS,n)$ complexes, and a
positive morphism $\overline{\theta}:\overline{H}\longrightarrow
\overline{\mathcal G}$ such that 
\[
(e,f,g)\in\overline{\theta}(\overline{H}_+)
\]
\end{lms}
\begin{demo}
First consider the case when $f=0$. We define $\overline{\theta}_g:\overline{H}_g\longrightarrow \overline{\mathcal G}$ by
\[
1\mapsto e; \overline{1}\mapsto \rho(e); 1\mapsto g
\]
on the $(\CS,n)$-complex, and note that
$(e,0,g)=\overline{\theta}(1,\overline{0},1)$. Similarly, when $g=0$,
we define  $\overline{\theta}_f:\overline{H}_f\longrightarrow \overline{\mathcal G}$
by
\[
1\mapsto e; (\overline{1},\overline{0})\mapsto \rho(e),
(\overline{0},\overline{1})\mapsto f; \overline{1}\mapsto \beta(f)
\]
on the $(\mathbb{I}^{\sim}_n,n)$-complex, and note that
$(e,f,0)=\overline{\theta}(1,(\overline{0},\overline{1}),\overline{0})$.

In the general case, we consider
\[
\xymatrix{
{\overline{H}_f\oplus \overline{H}_g}\ar[rr]^-{\overline{\theta}_f\oplus\overline{\theta}_g}\ar[d]_-{\overline{\gamma}}&&\overline{\mathcal G}\\
{\overline{H}}\ar[urr]_-{\overline{\lambda}}}
\]
with $\overline{\gamma},\overline{\lambda}$ chosen by Lemma
\ref{fiveone} above. By assumption,
\[
x=\gamma_0((1,0))=\gamma_0((0,1))
\]
so that
$(x,\gamma_n([(\overline{0},\overline{0}),(\overline{0},\overline{1})]),\gamma_1([1,\overline{0}]))$
is a positive preimage of $(e,f,g)$.
\end{demo}

\begin{thms}\label{alimit}
Fix $n\in\{2,3,\dots\}$, and let $\overline{\mathcal{G}}$ be an $n$-coefficient complex.  
Then, $\overline{\mathcal{G}} = \lim_{i \rightarrow \infty}
\overline{G}_i$, where each $\overline{G}_i$ is a finite direct sum of $(\mathbb{C},n)$,
$(\mathbb{I}^{\sim}_m,n)$ and $(\CS,n)$ complexes.
\end{thms}
\begin{demo}
Enumerate the positive elements of $\overline{\mathcal{G}}$ as
$(e_i,f_i,g_i)$ and apply Lemmas \ref{hit} and \ref{fiveone} alternately to get a diagram
\[
\xymatrix{
{\overline{H_1}}\ar[d]^-{\overline{\gamma}_1}\ar[rrr]^-{\overline{\theta}_1}&&&\overline{\mathcal{G}}\\
{\overline{G_1}}\ar@{^(->}[d]\ar[rrru]^-{\overline{\lambda}_1}&&&\\
{\overline{G_1}\oplus \overline{H_2}}\ar[d]^-{\overline{\gamma}_2}\ar[rrruu]^-{\overline{\lambda}_1\oplus\overline{\theta}_2}&&&\\
{\overline{G_2}}\ar@{^(->}[d]\ar[rrruuu]^-{\overline{\lambda}_2}&&&\\
{\vdots}&&&}
\]
where $(e_i,f_i,g_i)$ is the image under $\overline{\theta}_i$ of a positive element and 
$\ker{\overline{\theta}_i}=\ker{\overline{\gamma}_i}$. The maps will
then induce an order isomorphism.
\end{demo}

An inductive system of finite direct sums of building blocks is said
to have \emph{large denominators} if all connecting morphisms either
are zero on $\ko$ or have the $\kz$-component greater than or equal to $2$.

\begin{thms}\label{mainimho}
Let $n\in\{2,3,\dots\}$ and let $\overline{G}$ be a complex. The following are equivalent
\begin{enumerate}[(i)]
\item $\overline{G}$ is a $n$-coefficient complex;
\item $\overline{G}$ is an inductive limit of  finite direct sums of $(\mathbb{C},n)$,
$(\mathbb{I}^{\sim}_m,n)$ and $(\CS,n)$ complexes, and
$G_*$ has the Riesz property;
\item $\overline{G}$ is an inductive limit of  finite direct sums of $(\mathbb{C},n)$,
$(\mathbb{I}^{\sim}_m,n)$ and $(\CS,n)$ complexes, such
that the inductive system has
large denominators;
\item $\overline{G}\cong \BK(\A)$, where $\A$ is an $AD$
algebra of real rank zero. 
\end{enumerate}
\end{thms}
\begin{demo}
Note first that 
$(iv)$\infer $(i)$ was seen in Proposition \ref{isnzcc}.
Theorem \ref{alimit} proves $(i)$\infer $(ii)$, and since the property of
large denominators involves only the groups in $\gx$,
\cite[8.1]{gae:ccrrz} proves $(ii)$\infer$(iii)$. By compressing an
inductive system such as in $(iii)$ if necessary, we may assume that
each morphism among building blocks at level $i$ and level $i+1$ is
either zero on $\ko$ or greater than or
equal to $M_i$ on $\kz$, where $M_i$ is the largest number for which there is
an $(\DD{M_i},n)$ complex among the building blocks at level $i$. Then
by \cite{se:ciabtk} the inductive system can be realized by direct
sums of building blocks from the set
$\{\CC,\CS,\DD{2},\DD{3},\dots\}$ and  $*$-homomorphisms
among them. Furthermore, \cite[8.1]{gae:ccrrz} shows how to arrange
for real
rank zero in the limit.

By construction, the inductive limit $\A$ of this $C^*$-inductive
system is an AD algebra with the desired invariant
\[
\xymatrix{{\BK(\A):}&
{\K_0(\A)}\ar[r]&
{\K_0(\A;\ZZ/n)}\ar[r]&
{\K_1(\A).}
}
\]
However, since we have not used -- or even defined -- an ideal based
order on the building blocks $C(\mathrm{S}^1)$ and $\DD{m}$, cf.\
Definition \ref{egsemap}, we need to verify that the order on 
${\BK(\A)}$ coincides with the order on
$\overline{G}$. Since we have used the strict order on all the
algebraic building blocks this would follow directly if we knew that
all ideals of $\A$ arise as  inductive limits or direct sums of
subcollections of the building blocks in the system. And this in turn
is a consequence of the minimal real rank of $\A$, or directly by the
construction yielding this property in 
\cite[8.1]{gae:ccrrz}.
\end{demo}

It is essential to note at this stage that the ordered complex
$\BK(\A)$ is not complete for real rank zero AD algebras unless we
know that the torsion of $\ko$ is annihilated by the number $n$. Thus,
it is only in this case --- covered by \cite{se:ciabtk} --- that Theorem
\ref{mainimho} gives a one-to-one correspondence
between a class of \cstar-algebras and a class of algebraic invariants.

In this case, when the equivalent statements above hold true, we may
write the AD algebra as an inductive limit using only the buiding
blocks $C(\mathrm{S}^1)$ and $\mathbb{I}_n^\sim$, cf.\ \cite{se:ciabtk}

\section{The general case}
In the following we shall briefly recall definitions from \cite{mdse:ccpikc}.
Let $\Delta$ denote the ordered set $(\NN,\leq)$ where
\[
x\leq y\iff x\text{ divides }y.
\]
Note that $\Delta$ is directed, so that we may construct inductive
limits over $\Delta$. We will denote these by
\[
\Dilim{(G_p,f_{q,p})}
\]
where $f_{q,p}:G_p\rightarrow G_q$ are the bonding maps. 
When a cofinal subset $\Delta'$ of $\Delta$,
 is given, we may restrict attention to this, as
\[
\Dilim{G_n}\cong\ppDilim{G_n}.
\]

We define graded group homomorphisms
\[
\kmapk{mn}{m}:\pKu{m}{\A}\rightarrow\pKu{mn}{\A}
\]
by
\[
\left[\begin{matrix}\chi_{mn,n}&0\\0&\kmapk{mn}{m}\end{matrix}\right],
\]
where $\chi_{mn,n}$ is just multiplication by $m$ between the relevant
copies of $\K_0(\A)$. The maps $\kmapk{mn}{m}$ are positive, so we may define:

\begin{gather*}
\pKs{\A}=\Dilim\left({\pKu{n}{\A},\kmapk{mn}{n}}\right);\\
\pKs{\A}\po=\Dilim\left({\pKu{n}{\A}\po,\kmapk{mn}{n}}\right).
\end{gather*}

This gives the limit groups the structure of graded
ordered groups. The
even parts are naturally isomorphic to $\K_0(\A)\otimes \QQ$ since
\begin{eqnarray*}
\Dilim{(G,\comak{mn}{n})}&\cong&
\Dilim{(G\otimes\ZZ,\id\otimes\comak{mn}{n})}\\&\cong&
G\otimes\left({\Dilim(\ZZ,\comak{mn}{n})}\right)\cong
G\otimes \QQ
\end{eqnarray*}
naturally. We shall invoke this isomorphism tacitly in section 7.

The maps 
\[
\kmapk{mn}{m}:\pKu{m}{\A}\rightarrow\pKu{mn}{\A}
\]
can be described explicitly when $\A\in\{\CC,\CS,
\mathbb{I}_2^{\sim},\mathbb{I}_3^{\sim},\mathbb{I}_4^{\sim},\dots\}$
(\cite{se:iaa}).


In this section we study exact sequences of abelian groups
\begin{displaymath}
\GG\longrightarrow \GQ \stackrel{\rho}{\longrightarrow}
\mathcal{G}_{\mathbf{n}} \stackrel{\beta}{\longrightarrow} \mathcal{G}_1
\end{displaymath}
which are meant to represent the natural $\mathrm{K}$-theoretic invariants
for AD algebras of real rank zero having unbounded torsion in $\mathrm{K}_1$.
We begin by listing the properties that such an abstract sequence should
have before one may even consider whether the sequence arises as the invariant
\begin{equation}\label{Kninv}
\xymatrix{
{\BKf:}&&{\K_0(\A)}\ar[r]&
{\K_0(\A)\otimes\QQ}\ar[r]^-{{\rho}}&
{\K_0(\A;\QQ/\ZZ)}\ar[r]^-{{\beta}}&
{\K_1(\A)}}
\end{equation}
for some AD algebra $\A$ of real rank zero. We shall denote the
invariant consisting of two graded ordered groups and two group
homomorphisms as in \eqref{Kninv} by $\infinv(\A)$ and will see that
the conditions in the definition below are sufficient to ensure that the sequence
\begin{displaymath}
\GG \longrightarrow \mathcal{G}_0 \otimes \mathbb{Q} \stackrel{\rho}{\longrightarrow}
\mathcal{G}_{\mathbf{n}} \stackrel{\beta}{\longrightarrow} \mathcal{G}_1
\end{displaymath}
does indeed arise in such a manner.  

The next definition should be compared with Definition \ref{nccdef}.

\begin{defs}\label{gen-coeff-complex}
An exact sequence
\begin{displaymath}
\GG\longrightarrow\GQ \stackrel{\rho}{\longrightarrow}
\mathcal{G}_{\mathbf{n}} \stackrel{\beta}{\longrightarrow} \mathcal{G}_1
\end{displaymath}
(which we denote by $\overline{\mathcal{G}}$) of countably generated
abelian groups is an $\mathbf{n}$-coefficient complex if 
\begin{enumerate}[(i$\mbox{}_{\mathbf n}$)]
\item $\mathcal{G}_{\mathbf{n}}$ is pure torsion
\item $\im\beta = \tor \mathcal{G}_1$, and  every element $x \in \mathcal{G}_1$ of order $l$ has a $\beta$-lift of order $l$
\item $\mathcal{G}_* := \GG \oplus \mathcal{G}_1$ and $\mathcal{G}_\fnbox
:= \GQ \oplus \mathcal{G}_{\mathbf{n}}$ are graded ordered
groups inducing the same order on $\GG$
\item $\mathcal{G}_*$ has the Riesz interpolation property.
\item $\GQ \oplus \rho(\mathcal{G}_0)$ has the quotient order coming from 
$\mathrm{id}_{\mathcal{G}_0 \otimes \mathbb{Q}} \oplus \rho$
\item $\GG \oplus \beta(\mathcal{G}_{\mathbf{n}})$ has the quotient order coming from 
$\mathrm{id}_{\mathcal{G}_0} \oplus \beta$ 
\item $\mathcal{G}_0$ is unperforated and $\mathcal{G}_*$ is weakly unperforated.
\end{enumerate}
\end{defs}

Note that since $\GG$ is torsion free, it is determined by $\GQ$ in
this setup. In the proofs below we may hence concentrate our work on
the rightmost three groups in the complex.

\begin{props}\label{decomposition}
Let
\begin{displaymath}
\overline{\mathcal{G}}:\qquad \GQ \stackrel{\rho}{\longrightarrow}
\mathcal{G}_{\mathbf{n}} \stackrel{\beta}{\longrightarrow} \mathcal{G}_1
\end{displaymath}
be an $\mathbf{n}$-coefficient complex.  Then, there exist a strictly increasing
sequence of natural numbers $\langle n_i \rangle$ with $n_i|n_{i+1}$, a
$n_i$-coefficient complex $\overline{G^i}$ for each $i$, and positive morphisms
$\overline{\theta_i}:\overline{G^i} \rightarrow \overline{G^{i+1}}$ such that
\[
\overline{\mathcal{G}} \simeq (\overline{G^i},\overline{\theta_i}).
\]
\end{props}

\begin{demo}
Let $\langle n_i \rangle$ be a strictly increasing sequence of natural 
numbers with the property that every natural number divides some $n_i$.
By the main theorem of section 5 of \cite{gae:dgt}, the graded ordered group
$(\GQ,\mathcal{G}_1)$ in Definition \ref{gen-coeff-complex} is the limit of an 
inductive sequence
\[
\left((G_0^i,G_1^i),(\phi_0^i,\phi_1^i) \right)
\]
of graded ordered groups $(G_0^i,G_1^i)$.  Furthermore, each $(G_0^i,G_1^i)$
 consists of the first and
third groups of an $n_i$-coefficient complex which is a direct sum
of $(\mathbb{C},n_i)$, $(\mathbb{I}^{\sim}_m,n_i)$ and $(C(\mathrm{S}^1),n_i)$
complexes.  Let 
\[
\phi_0^{i \infty}:G_0^i \rightarrow \mathcal{G}_0, \ \ \phi_1^{i \infty}:G_1^i 
\rightarrow \mathcal{G}_1,
\]
be the canonical maps.  Assume that the $n_i$ have been chosen large enough for 
the inclusions $G_0^i \subseteq \rho^{-1}(\mathcal{G}_{\mathbf{n}}[n_i])$ and
$\tor(G_1^i) \subseteq \beta(\mathcal{G}_{\mathbf{n}}[n_i])$ to hold.  
We then have
\[
(\GQ,\mathcal{G}_1) = \lim_{i \rightarrow \infty}
\left((\rho^{-1}(\mathcal{G}_{\mathbf{n}}[n_i]),\beta(\mathcal{G}_{\mathbf{n}}[n_i]) 
\cup \phi_1^{i \infty}(G_1^i)),\iota \right),
\]
where $\iota$ is the inclusion map.
Note that
\[
(\rho^{-1}(\mathcal{G}_{\mathbf{n}}[n_i]),\beta(\mathcal{G}_{\mathbf{n}}[n_i]) 
\cup \phi_1^{i \infty}(G_1^i))
\]
is a graded ordered group for each natural number $i$ --- it is an order hereditary
subgroup of $(\GQ,\mathcal{G}_1)$.  One can then verify that the complex
\[
\rho^{-1}(\mathcal{G}_{\mathbf{n}}[n_i]) \stackrel{\rho}{\rightarrow}
\mathcal{G}_{\mathbf{n}}[n_i] 
\stackrel{\beta}{\rightarrow} \beta(\mathcal{G}_{\mathbf{n}}[n_i]) \cup \phi_1^{i \infty}(G_1^i) 
\]
is an $n_i$-coefficient complex.  (The only subtle point is the exactness of the sequence, 
which follows from the second half of property $iii_{\mathbf n}$ of Definition \ref{gen-coeff-complex}.)  
The limit of the inductive sequence
\[
\left(\rho^{-1}(\mathcal{G}_{\mathbf{n}}[n_i]) \stackrel{\rho}{\rightarrow}
\mathcal{G}_{\mathbf{n}}[n_i] \stackrel{\beta}{\rightarrow} 
\beta(\mathcal{G}_{\mathbf{n}}[n_i]) \cup \phi_1^{i \infty}(G_1^i), \kappa_{i,i+1} \right)
\]
--- $\kappa_{i,i+1}$ is the inclusion map --- is then 
\[
\GQ \stackrel{\rho}{\longrightarrow}
\mathcal{G}_{\mathbf{n}} \stackrel{\beta}{\longrightarrow} \mathcal{G}_1
\] 
by construction.
\end{demo}


The following theorem is the generalised integer coefficient version of
Theorem 5.3.

\begin{thms}\label{seventwo}
Let $\overline{\mathcal{G}}$ be an $\mathbf{n}$-coefficient complex.  
Then, 
\[
\overline{\mathcal{G}} = \lim_{i \rightarrow \infty}
(\overline{H}_i,\overline{\gamma_i}), 
\]
where each $\overline{H}_i$ is a direct sum of $(\mathbb{C},n_i)$,
$(\mathbb{I}^{\sim}_m,n_i)$ and $(\CS,n_i)$ complexes, some $n_i \in \mathbb{N}$.  
\end{thms}

\begin{demo}
Assume the inductive sequence decomposition of $\overline{\mathcal{G}}$ from Lemma 7.2.  
For brevity, write
\[
\overline{H}^{n_i} = \rho^{-1}(\mathcal{G}_{\mathbf{n}}[n_i]) \stackrel{\rho}{\rightarrow}
\mathcal{G}_{\mathbf{n}}[n_i] 
\stackrel{\beta}{\rightarrow} \beta(\mathcal{G}_{\mathbf{n}}[n_i]) \cup \phi_1^{i \infty}(G_1^i) 
\]
  Put $\kappa_{l,m} = \kappa_{m-1,m} \circ
\cdots \circ \kappa_{l,l+1}$.
By Theorem \ref{alimit}, each $\overline{H}^{n_i}$ is the limit of an inductive
system $(\overline{H}^{n_i}_k,\overline{\theta}^{n_i}_{k,k+1})$, where each
$\overline{H}^{n_i}_k$ is a direct sum of $(\mathbb{C},n_i)$,
$(\mathrm{I}^{\sim}_m,n_i)$ and $(\mathrm{S}^1,n_i)$ complexes.  

It will suffice to define a sequence of positive morphisms
\[
\gamma_{i,i+1}:\overline{H}_i^{n_i} \rightarrow \overline{H}_{i+1}^{n_{i+1}}
\] 
making the diagram
\[
\xymatrix{
&&&&&&  \overline{\mathcal{G}} \\
&&&&&&\\
&&&&&& \mbox{}\ar@{.}[uu] \\
&&&&&&\\
{\overline{H}_1^{n_2}\ar[rr]^{\overline{\theta}_{1,2}^{n_2}}} &&
{\overline{H}_2^{n_2}\ar[rr]^{\overline{\theta}_{2,3}^{n_2}}\ar[uurr]^{\gamma_{2,3}}} && \mbox{}\ar@{.}[rr] && 
{\overline{H}^{n_2} \ar[uu]^{\kappa_{2,3}}}\\
&&&&&&\\
{\overline{H}_1^{n_1}\ar[rr]^{\overline{\theta}_{1,2}^{n_1}}}\ar[uurr]^{\gamma_{1,2}} &&
{\overline{H}_2^{n_1}\ar[rr]^{\overline{\theta}_{2,3}^{n_1}}} && \mbox{}\ar@{.}[rr] && 
{{\overline{H}^{n_1}}\ar[uu]^{\kappa_{1,2}}}
}
\]
commute;  by compressing the sequence for $H^{n_i}$ one can ensure that every positive 
element in $\mathcal{G}$ is the image under $\gamma_{i,\infty}$ of a positive element 
in $\overline{H}_i^{n_i}$ for some $i \in \mathbb{N}$.  

Let $\overline{\theta}_{1,\infty}^{n_1}:\overline{H}_1^{n_1} \rightarrow \overline{H}^{n_1}$
be the canonical morphism.  Let $M$ be a minimal set of positive generators for 
$\overline{H}_1^{n_1}$.  Find, 
by compressing the inductive sequence for $\overline{H}^{n_2}$ if necessary, a set
$\widetilde{M}$ of positive pre-images via $\overline{\theta}_{1,\infty}^{n_2}$ of the 
elements of $\kappa_{1,2} \circ \overline{\theta_{1,\infty}}(M)$ in $\overline{H}_2^{n_2}$.
Note for future reference that each element of $\kappa_{1,2} \circ \overline{\theta_{1,\infty}}(M)$ is divisible by
$n_{i+1}/n_i$ inside $\overline{H}^{n_2}$, so we may assume that $m(n_i/n_{i+1}) \in
\overline{H}_2^{n_2}$ whenever $m \in \widetilde{M}$. 
Define $\gamma_{1,2}$ by sending an element $m \in M$ to the corresponding pre-image
of $\kappa_{1,2} \circ \overline{\theta_{1,\infty}}(m)$ in $\widetilde{M}$.

\end{demo}

The following theorem is the generalised integer coefficients
version of Theorem \ref{mainimho}. 

\begin{thms}\label{seventhree}
The following are equivalent:
\begin{enumerate}[(i)]
\item $\overline{\mathcal{G}}$ is a $\mathbf{n}$-coefficient complex;
\item $\overline{\mathcal{G}}$ is an inductive limit of  finite direct sums of $(\mathbb{C},n)$,
$(\mathbb{I}^{\sim}_m,n)$ and $(C(\mathrm{S}^1),n)$ complexes, where $n$ ranges over the 
natural numbers, and $\mathcal{G}_*$ is a Riesz group;
\item $\overline{\mathcal{G}}$ is an inductive limit of  finite direct sums of $(\mathbb{C},n)$,
$(\mathbb{I}^{\sim}_m,n)$ and $(C(\mathrm{S}^1),n)$ complexes where $n$ ranges over the
natural numbers, and such that the inductive system has large denominators;
\item $\overline{\mathcal{G}}\cong \infinv(\A)$, where $\A$ is an $AD$
algebra of real rank zero. 
\end{enumerate}
\end{thms} 
\begin{demo}
The proof follows the proof of Theorem \ref{mainimho} with the
exception that in $(iii)\infer(iv)$, one needs to realize maps from
$\BK[n_i](\A_i)$ to $\BK[n_{i+1}](\A_i)$ by a triple of maps of the
form
\[
(\chi_{n_{i+1},n_i}\circ f_*,\kappa_{n_{i+1},n_i}\circ f_*,f_*)
\]
rather than directly by a $*$-homomorphism. However, as noted at the
end of the proof of Theorem \ref{seventhree}, the maps in question
will have a $\K_0$-component which is divisible by $n_{i+1}/n_i$, so
this may be arranged as in the proof of Theorem \ref{mainimho}.
\end{demo}

\section{The example of Dadarlat and Loring}
It follows from the work of B\"odigheimer
(\cite{cfb:skskI},\cite{cfb:skskII}) that 
the unspliced short exact sequence
\[
\xymatrix{{0}\ar[r]&
{\K_0(\A)/n}\ar[r]&
{\K_0(\A;\ZZ/n)}\ar[r]&
{\K_1(\A)[n]}\ar[r]&{0}
}
\]
will always split. This has been useful in the analysis of other
aspects of this object (\cite{se:ksccfmi}, \cite{se:ahcfmi}) but we
have not been able to employ the fact in the proofs leading to Theorem
\ref{mainimho} -- since the splitting map is unnatural, it is
difficult to use it when trying to establish the range of the invariant.
By contrast, it is a useful result when trying to describe the amount
of freedom one has in the choice of equipping $\overline{G}$ as an
$n$-coefficient complex when $\gx$ is fixed, as we shall se below.

However, since Theorem
\ref{mainimho} combines with the results mentioned above to prove that
every $n$-coefficient complex will split when unspliced, it is perhaps worthwhile to
note that this follows already from properties $(i)$, $(ii)$ and the fact
(contained in $(vii)$) that $G_0$ is torsion free. One proves this by
first establishing that $\im\rho$ is a pure subgroup of $\gd$ and then
appealing to \cite{lf:iagI}.

\begin{remar}
For $G_0$ not necessarily torsion
free the properties $(i)$--(vi), $(viii)$-(ix) would not imply splitness, as the example
\[
\xymatrix{
{\ZZ[\frac14]\oplus\ZZ/2}\ar[r]&
{\ZZ/4}\ar[r]&{\ZZ/2}}
\]
shows. When this is equipped with the strict order induced by the
standard order on $\ZZ[\frac14]$ it has all the properties of our $4$-coefficient complexes
except unperforation, but could not be the augmented $\K$-theory of a
\cstar-algebra. Thus to extend range results beyond the case
considered above, on would have to impose an extra condition; for
instance that $\im\rho$ was a pure subgroup of $\gn$.
\end{remar}

Using the splitness, we may, up to isomorphism, for any $n$-coefficient complex write $\gd=R\oplus B$ such
that $R=\gz/n$, $B=\go[n]$ and 
\[
\rho(x)=(x+n\gz,0)\qquad
\beta(r,b)=b
\]
The properties $(i)$-(ix) simplify accordingly.

\begin{remar}
Note that graded ordered every dimension group with torsion $\gx$ can
be extended to an $n$-coefficient complex when $G_0$ is unperforated, simply by ordering
$G_0\oplus(R\oplus B)$ by
\[
(x,(r,b))\Longleftrightarrow\begin{cases}(x,(r,0))\geq 0\\(x,b)\geq
  0\end{cases}
\]
where the quotient order of $G_0\oplus G_0$ and the order on $\gx$,
respectively, are used to determine
whether $(x,(r,0))$ and $(x,b)$ are positive. 
\end{remar}

Fix $\gx$. We have seen that up to isomorphism, every $n$-coefficient complex is of the form
\[
\xymatrix{
{\gz}\ar[r]&
{R\oplus B}\ar[r]&
{\go,}
}
\]
so determining how many \nzccs with this particular $\gx$ are possible
comes out to determining which order structures on $\gd$ will satisfy
properties $(iii)$, $(v)$, $(vi)$, $(viii)$ and $(ix)$. 

As an example this process, and an application of Theorem \ref{mainimho}, let us
return to the example of Dadarlat and Loring which originally established the
need for ordered $\K$-theory with coefficients. They considered the
$\gx$-group given by
\begin{gather}
\gz=\left\{(x,y_i)\in\ZZ[\tfrac1{n+1}]\oplus\prod_{-\infty}^{\infty}\ZZ\mid
  y_i=x(n+1)^{|i|}, a.e. (i)\right\}\label{dl}\\
\go=\ZZ/n\notag
\end{gather}
equipped with the standard order on $\gz$ and the strict order
herefrom on $\gx$. In \cite{mdtal:ccomk} examples were given to show,
in effect,
that there were two different ways to complete $\gx$ to an $n$-coefficient complex.

For convenience, let $n=2$. With the notation above we have
\begin{gather*}
R=\left\{(a,b_i)\in\ZZ/2\oplus\prod_{-\infty}^{\infty}\ZZ/2\mid
  b_i=a, a.e. (i)\right\}\\
B=\{c\in\ZZ/2\}
\end{gather*}
where
\[
\rho(x,y_i)=(\overline{a},\overline{y_i})\text{ with }x=\frac a{3^i}.
\]

In the proof below we need elements $\delta_j,\Delta_N\in(\ZZ/2)^\ZZ$
defined by
\[
(\delta_j)_i=\begin{cases}1&i=j\\0&i\not=j\end{cases}\qquad
(\Delta_N)_i=\begin{cases}1&|i|> N\\0&|i|\leq N\end{cases}
\]

\begin{props} Up to isomorphism every \nzcc[2]
completing $\gx$ defined in \eqref{dl} will order $\gd$ by
\[
((x,y_i),(a,b_i,c))\geq 0\Longleftrightarrow
\begin{cases} x\vdash a\\y_i\vdash b_i+\epsilon_ic\\x\vdash c\end{cases}
\]
where $(\epsilon_i)_{i\in\ZZ}$ is a sequence
in $(\ZZ/2)^\ZZ$.
\end{props}
\begin{demo}
It is easy to see that this defines a \nzcc[2]. In the other direction,
first note that since the order of $\gz\oplus (R\oplus 0)$ is induced by the graded
ordered group $\gz\oplus\gz$ we have
\[
((x,y_i),(a,b_i,0))\geq 0\Longleftrightarrow
\begin{cases} x\vdash a\\y_i\vdash b_i\end{cases}.
\]

Lift the element $((1,3^{|i|}(1-\delta_j)),1)\in\gx$ to a positive
element $((1,3^{|i|}(1-\delta_j)),(a,b_i,1)\in\gd$. If $a=1$, note
that also
\[
((1,3^{|i|}(1-\delta_j)),(a+1,b_i+1-\delta_j,1+0))
\]
is positive, so that we may without loss of generality assume that
$a=0$. Similarly, we may assume that $b_i=0$ for all $i\not=j$.

We have hence seen that at least  one of 
\[
((1,3^{|i|}(1-\delta_j)),(0,\delta_j,1))\qquad
((1,3^{|i|}(1-\delta_j)),(0,0,1))\qquad
\]
is positive in $\gd$. We will define $\epsilon_j$
accordingly such that
\[
((1,3^{|i|}(1-\delta_j)),(0,\epsilon_j\delta_j,1))\geq 0.
\]

If $((x,y_i),(a,b_i,1))\geq 0$ then $x>0$ because $(x,1)\geq0$ in
$\gx$. Further, if $y_j=0$ then since
\[
((1+x,3^{|i|}(1-\delta_j)+y_i),(a,b_i+\epsilon_i\delta_j,1+1))\geq 0 
\]
we have that $b_j+\epsilon_j=0$. We conclude that $x\vdash a,c$ and
$y_j\vdash b_j+\epsilon_jc_j$.

In the other direction, assume that $x\vdash a,c$ and $y_j\vdash
b_j+\epsilon_jc$. We already know that $((x,y_i),(a,b_i,c))\geq 0$
when $c=0$, so we can focus on the case $c=1$. In this case we will
have $x>0$ and hence $y_i>0$ for $|j|\geq N$ for some $N\in\NN$. We lift $((x,y_i),1)$
to some positive element $((x,y_i),(a',b_i',1)$. If $a\not=a'$, we note that
$((x,y_i),(1,\Delta_N,0))\geq 0$ 
whence also
\[
((x,y_i),(a'+1,b_i'+\Delta_N,1))\geq 0
\]
so that we without loss of generality may assume that $a=a'$ and hence that
$b_i=b_i'$ for all but finitely many $i$. For the remaining $i$s, if
$y_i>0$, we can adjust to get $b_j=b_j'$ in a similar fashion. And if
$y_i=0$ then since we have
\[
((x,y_i)(a+a',b_j+b_j',1+1))=((x,y_i)(0,b_j+b_j',0))\geq 0
\]
we get that $b_j=b_j'$.
\end{demo}

Let ``$\sim$'' be the finest equivalence relation on $(\ZZ/2)^\ZZ$ such
that
\[
(\epsilon_j)\sim(\epsilon_{j+k})\sim(1+\epsilon_{j})\sim(\epsilon_j\Delta_N(j))\sim(\epsilon_{(-1)^{r_{|j|}}j})
\]
for any $k\in\ZZ$, $N\in\NN$ and any
$r_j\in(\ZZ/2)^{\NN\cup\{0\}}$. Using methods from \cite{mdtal:ccomk}
one can prove that the augmentations associated to $(\epsilon_i)$ and
$(\eta_i)$ in $(\ZZ/2)^\ZZ$ are isomorphic precisely when
$(\epsilon_i)\sim(\eta_i)$. The examples given in \cite{mdtal:ccomk}
correspond to $(\epsilon_i)=(0)$ and $(\eta_i)$ given by $1$ on
positive entries and $0$ on negative ones.

One sees easily that there are uncountably many nonisomorphic
2-coefficient complexes in this case --- even though we have only added one bit of
information to $\gz\oplus(R\oplus 0)$ the amount of freedom in
choosing the order structure is immense.

\providecommand{\bysame}{\leavevmode\hbox to3em{\hrulefill}\thinspace}
\providecommand{\MR}{\relax\ifhmode\unskip\space\fi MR }
\providecommand{\MRhref}[2]{%
  \href{http://www.ams.org/mathscinet-getitem?mr=#1}{#2}
}
\providecommand{\href}[2]{#2}

\end{document}